\documentclass[12pt]{article}

\title{The Stationary Measure of a 2-type Totally Asymmetric Exclusion
Process}
\author{Omer Angel\thanks{Laboratoire de Math., Universite Paris-Sud,
Bat. 425, 91405 Orsay Cedex, France. {\tt angel@math.ubc.ca}}}
\date{July 2004}

\usepackage{amsmath, amsthm, amssymb,amsfonts,mathbbol}

\newtheorem{thm}{Theorem}[section]
\newtheorem{lemma}[thm]{Lemma}
\newtheorem{coro}[thm]{Corollary}
\newtheorem{conj}{Conjecture}
\theoremstyle{definition}
\newtheorem{defn}[thm]{Definition} 

\newcommand{\lemref}[1]{Lemma~\ref{lem:#1}}
\newcommand{\thmref}[1]{Theorem~\ref{thm:#1}}
\newcommand{\corref}[1]{Corollary~\ref{cor:#1}}

\newcommand{\Z}{\ensuremath{\mathbb{Z}}}
\newcommand{\E}{\ensuremath{\mathbb{E}}}
\renewcommand{\P}{\ensuremath{\mathbb{P}}}

\newcommand{\dom}{\succcurlyeq}
\newcommand{\seq}[1]{{\text{\tt #1}}}
\newcommand{\zero}{\seq{0}}
\newcommand{\one}{\seq{1}}
\renewcommand{\star}{\seq{*}}

\renewcommand{\chi}{\mathbb{1}}
\begin{document}

\maketitle
\begin{abstract}
  We give a combinatorial description of the stationary measure for a
  totally asymmetric exclusion process (TASEP) with second class particles,
  on either $\Z$ or on the cycle $\Z_N$. The measure is the image by a
  simple operation of the uniform measure on some larger finite state
  space. This reveals a combinatorial structure at work behind several
  results on the TASEP with second class particles.
\end{abstract}

\section{Introduction}
%%%%%%%%%%%%%%%%%%%%%%

In the normal totally asymmetric simple exclusion process (TASEP), a number
of particles occupy some distinct vertices of a graph, which in our case
will be taken to be either $\Z$ or the cycle $\Z_N$. Each particle moves to
the position to its right at rate one whenever that position is unoccupied.
This gives a Markov chain with continuous time on the set of configurations
--- sets of locations for the particles. It is well known (and easy to
prove) that if there are $a$ particles on $\Z_N$, then the stationary
distribution is uniform among all $N\choose a$ possible states. Since in
the cyclic case the total number of particles is invariant, any other
stationary measure is some linear combination of these uniform measures.
Since the set of stationary measures is convex, from now on we focus only
on the extremal stationary measures.

In the infinite case of $\Z$ there are two types of extremal stationary
measures. There are the so called trivial stationary measures where the
particles occupy all positions to the right of some point and no other
positions. The non-trivial stationary measures are those where each vertex
is occupied with some probability $p\in[0,1]$ independently of all other
vertices. Note that these measures are limits as $N\to\infty$ of the unique
stationary measure on $\Z_N$ with $[pN]$ particles. For additional
background on the TASEP and many references see \cite{liggett1,liggett2}.

We are interested in the process with second class particles defined as
follows (see e.g. \cite{DJLS}). The particles are classified as first class
and second class particles. Sites are in one of three states: empty, or
occupied by a single particle of one class or the other. As with the simple
process, each particle jumps to the position to its right at rate 1 when
that position is empty. Additionally, whenever a first class particle has a
second class particle to its right, the two swap places at rate 1 (thus
second class particles may move in both directions).

\medskip

It is interesting at this point to note that the TASEP on a cycle is
equivalent to the following shuffling method of cards. The cards are
arranged in a cycle. Each consecutive pair is chosen at rate 1, and the
pair is sorted with the larger card to the right. Assume there are $a$ high
cards, $b$ medium cards and $c$ low cards, and that cards within the same
class are not distinguished from one another. The dynamics are exactly
equivalent to those of the exclusion process with high cards corresponding to first
class particles, medium cards to second class particles and low cards to
empty positions (see \cite{BBHM} for an application). It is of course
interesting to study the case where there are more than three types of
cards (or particles). These cases seem significantly harder to analyze,
though experimental studies of small cases, as well as the extremal case of
$N$ different particles do show interesting phenomena.

We give one final equivalent interpretation of the process, that may be
related to the fact that adding further classes of particles breaks known
techniques for working with the model. Non-empty sites of the graph are
occupied by either a particle or an anti-particle. Each {\bf edge} is
chosen at rate 1. If a particle can move right across the edge to an empty
spot it does. If an anti-particle can move left across the edge to an empty
space it does. Finally, if there is a particle on the left and an
anti-particle on the right, they both move, exchanging their positions, and
otherwise nothing happens. Thus there are particles moving right and
anti-particles moving left, but the movements have rate 1 at each edge
rather than each particle. This is equivalent to the previous form, with
anti-particles representing empty spaces, and empty spaces representing the
second class particles. When writing out states of the process, we will use
\seq{1}'s for particles, \seq{0}'s for anti-particles and \seq{*}'s for
empty spaces. We will mostly refer to this interpretation of the process
from now on, and use the term TASEP to refer to this process as well.

If the initial state is such that only two of the three types of positions
appear (i.e.\ if there are no particles of one type or the other, or if all
sites are occupied), then the process is simply the well understood
exclusion process. However, the relation to the 1-type case is deeper than
that. If an observer sees only the particles, and ignores the
anti-particles, treating them as empty spaces, then he observes a regular
exclusion process. Thus in the stationary distribution with $a$ particles
and $b$ anti-particles, the marginal of the positions of particles is
uniform over all $N\choose a$ sets. Similarly, the anti-particles form an
exclusion process on their own (up to a reversal of the directions) and are
therefore uniformly distributed over $N\choose b$ possible sets. Of course
the two sets of locations are disjoint and hence dependent.

A similar statement holds in the infinite space case. In a stationary
measure, if one type of particles are in a trivial state, then so must be
the other, and we get the trivial stationary measures where a single
segment of empty spaces has particles on its right and anti-particles on
its left. Otherwise, the set of positions of particles have an i.i.d.\
product measure with some density, as do the anti-particles. Of course, the
two marginals are not independent, nor is their joint distribution a
product measure.

The main result of this paper is a combinatorial description of the
stationary measures for the TASEP with three types of positions (empty and
two types of particles). These results are then used to shed light on some
known properties of the stationary measures. Duchi and Schaeffer \cite{DS}
have found a similar relation, and use it to derive other results on the
TASEP (primarily on a finite interval).

\section{Definitions and Results}
%%%%%%%%%%%%%%%%%%%%%%%%%%%%%%%%%

\begin{defn} \label{def:collapse}
Two sets of positions $S,T$, of positions on $\Z$ or $\Z_N$, (not
necessarily disjoint), are said to {\em collapse} to a state $x$ of the
TASEP if $x$ is the result of the following (collapsing) procedure:
Anti-particles are placed at the locations specified by $T$. Next, the
locations in $S$ are checked (in an arbitrary order). If a location is
empty, a particle is placed there. Otherwise a particle is placed in the
nearest empty position {\bf to the left} of the specified location.
\end{defn}

In the case of the cycle $\Z_N$, we will only use the process for sets
$S,T$ with $|S|+|T|\le n$, so that there is always an open position for
every particle in $S$. In the case of $\Z$, if there are no open positions
to the left of some element in $S$, we disregard that element.

The collapsing procedure defines a function from the product space (pairs
of sets) to the space of states of the TASEP. The order at which the
positions of $S$ are used has no effect on the final resulting state $x$:
The anti-particles in $x$ are located exactly in the positions in $T$. A
position $a$ contains a particle in the resulting state $x$ if and only if
$a \notin T$ and there is some interval $I=[a,b]$ such that $|I\cap S| +
|I\cap T| \ge |I|$. These statements may be used as an equivalent condition
of collapsing, avoiding the possibly infinite algorithmic definition in the
case of $\Z$.

\begin{thm} \label{thm:cycle}
The stationary measure for the exclusion process on $\Z_N$ with $a$
particles and $b$ anti-particles is the image by collapsing of the uniform
measure on pairs of sub-sets $S,T$ of the cycle of sizes $a,b$ respectively.
\end{thm}

As a corollary, it follows that the least likely states for given
population sizes are \seq{1\ldots1*\ldots*0\ldots0} and its cyclic shifts,
i.e.\ positions where sites of each type form a single interval, with
positive particles to the left of the empty sites (and anti-particles to
the right). These states have probability $1/{N \choose a}{N\choose b}$
each, while every other state has a probability that is some integer
multiple of this probability. Section~\ref{sec:open} contains some
conjectures generalizing these facts to the case of more particle types.

As noted above, On $\Z$ the TASEP with second class particles has trivial
stationary measures, which correspond to the states with particles to right
of some point, and anti-particles to the left of some other point with
empty spaces in between. In the notation of first and second class
particles these are states with second class particles in some interval,
first class particles to the right and empty spaces to the left.

\begin{thm} \label{thm:line}
The non-trivial extremal invariant measures for the exclusion process on
$\Z$ are the image by collapsing of sets $S,T$ where each $n\in S$ with
probability $p$ and $n\in T$ with probability $q$, all independently.
\end{thm}

Note that if $p+q\ge1$, then a.s.\ all positions of $\Z$ will eventually
have some particle on them, and then the resulting measure is one where
each site has an particle with probability $1-q$, and an anti-particle
otherwise, independently of all other sites. As noted above, in the lack of
empty spaces these are the (non-trivial extremal) invariant measures. If
$p+q<1$, then the resulting measure has particles with density $p$,
anti-particles with density $q$ and empty sites with density $1-p-q$.

Theorems~\ref{thm:cycle} and \ref{thm:line} can also be used to derive
properties of the collapsing procedure. Since the distribution of the
positions of particles in the stationary distribution is known:

\begin{coro}
Let $S'$ is the set of locations of particles after collapsing independent
$S,T$, which are have either the uniform distribution over sets of a given
size (in $\Z_N$) or i.i.d.\ product distributions (on $\Z$). Then $S'$ also
has a uniform or an i.i.d.\ product distribution.
\end{coro}

Furthermore, the TASEP itself has a natural symmetry of reversing the
charge of the particles as well as the direction of the graph. It follows
that a similar ``dual'' collapsing procedure where the particles are fixed
and the anti-particles are moved forward to empty spaces also yields the
same stationary measure.

\medskip

In the next section we prove a combinatorial lemma that is closely related
to stationarity of the collapsed uniform measures. Sections~\ref{sec:cycle}
and \ref{sec:line} contains the proof of Theorems~\ref{thm:cycle}, and
\ref{thm:line} respectively. In Section~\ref{sec:line} we also use the
collapsing description of the stationary measures to shed new light on some
of the results of \cite{FFK}. Finally, Section~\ref{sec:open} contains some
open problems and conjectures regarding more general multi-type asymmetric
exclusion processes.

\section{Binary trees and dominated sequences}
%%%%%%%%%%%%%%%%%%%%%%%%%%%%%%%%%%%%%%%%%%%%%%

This section contains the combinatorial foundation for proving the
stationarity of the collapsed uniform measure. The key result here is a
bijection between binary trees and {\bf pairs} of binary sequences, (i.e.\
sequences made up of \zero's and \one's). The sequences are related to the
TASEP since a binary sequence describes a segment in a state of the TASEP
with no empty sites (which is one of the reasons we use the
$\{\seq{0,*,1}\}$ notation).

\begin{defn}
Consider two finite binary sequences $A,B$ of the same length $n$. We say
that $A$ {\em dominates} $B$ and write $A\dom B$ if it is possible to get
from $A$ to $B$ by moving \one's to the right. The {\em weight} of a binary
sequence $A$ is defined as the number of binary sequences dominated by it:
\[
W(A) = | \{ B: A\dom B \} | .
\]
\end{defn}

In particular, it is necessary for $A\dom B$ that both sequences have the
same number of ones. Let the number of \one's in the first $i$ digits of
$A$ (resp. $B$) be denoted by $a(i)$ (resp. $b(i)$). A condition equivalent
to $A\dom B$ is that $a(i)\ge b(i)$ for all $i$, and $a(n)=b(n)$. Thus for
example, $W(\seq{1010})=5$. If $x=\seq{1\ldots10\ldots0}$ has $k$ ones
followed by $l$ zeroes, then $W(x)= {k+l \choose k}$. If the \seq{1}'s and
\seq{0}'s were exchanged $W(x)$ would have been 1.

We use the following version of binary trees:

\begin{defn}
A {\em binary tree} is a rooted tree where each vertex including the root
may have a left child, marked as left, and may have a right child, marked
as right. A vertex may have either child, both, or neither.
\end{defn}

Thus a tree is either the empty tree with only a root vertex, or it has a
left sub-tree, a right sub-tree or both. Note that having an empty sub-tree
on some side is different from not having a sub-tree on that side.

Next, we define recursively a function $f$ mapping binary trees to binary
sequences, as follows. The empty tree (with no edges) encodes the empty
sequence. Otherwise,
\[
f(T) = \left\{ \begin{array}{ll}
        f(L)\zero          & \text{if $T$ has only a left sub-tree $L$,} \\
        \one f(R)          & \text{if $T$ has only a right sub-tree $R$,} \\
        f(L)\seq{01}f(R)   & \text{if $T$ has sub-trees $L,R$,} \\
\end{array} \right.
\]
where e.g.\ $f(L)\seq{01}f(R)$ means a concatenation of $f(L)$, \seq{01}
and $f(R)$.
In this way any binary 
sequence may be encoded by a binary tree, though generally the encoding is
not unique. The length of the sequence is the number of edges of the tree,
and the number of \one's is the total number of right children. The
following combinatorial Lemma and immediate Corollary relate binary trees
and dominated sequences:

\begin{lemma} \label{lem:LR_bij}
There exists a bijection between binary trees $T$ and pairs $A,B$ of binary
sequences such that $A\dom B$ and $A=f(T)$.
\end{lemma}

\begin{coro} \label{cor:LR_w}
The number of binary trees that encode a given binary sequence $A$ is
$W(A)$.
\end{coro}

\begin{proof}[Proof of Lemma~\lemref{LR_bij}]
First, we define in a similar fashion a second function $g$, that maps
binary trees to binary sequences. The empty tree is mapped to the empty
sequence, and otherwise, 
\[
g(T) = \left\{ \begin{array}{ll}
        \zero g(L)          & \text{if $T$ has only a left sub-tree $L$,} \\
        \one g(R)           & \text{if $T$ has only a right sub-tree $R$,} \\
        \zero g(L)\one g(R)   & \text{if $T$ has sub-trees $L,R$,} \\
\end{array} \right.
\]
We now show that mapping a tree $T$ to a pair of sequences $A=f(T)$,
$B=g(T)$ gives a bijection satisfying the Lemma's requirements.

The following facts are clear: for sequences $A,B,C,D$ such that $A\dom B$
and $C\dom D$ we have $A\zero \dom \zero B$ and $AC \dom BD$. By induction,
it follows that $f(T)\dom g(T)$: If the root has only a right child, then
$f(T) = \one f(R) \dom \one g(R) = g(T)$. If the root has only a left
child, then $f(T) = f(L) \zero \dom \zero g(L) = g(T)$. If the root has two
offspring, then concatenation of the previous two cases yields $f(T) \dom
g(T)$.

To see that this mapping is a bijection, we show how to (recursively)
recover from sequences $A,B$ a tree that is mapped to them, in such a way
that at each step there is a unique possibility, so that the tree is
unique. 

Consider a binary tree $T$, with two sub-trees $L,R$ and sequences
$A=f(T)=f(L)\seq{01}f(R)$ and $B=g(T)=\zero g(L)\one g(R)$. Recall that $a(i)$
and $b(i)$ count \one's in preambles of $A$ and $B$. Since $f(L)\dom g(L)$,
we have that $a(i)\ge b(i+1)$ for $i\le|L|$. However, $a(|L|+1)=b(|L|+1)$
is the number of right edges in $L$, and the next bit in $B$ is a one. Thus
the inequality $a(i)\ge b(i+1)$ fails for for the first time for $i=|L|+1$.
The following algorithm emerges: Given $A,B$ find the first $i$ for which
$a(i)<b(i+1)$, and set $|L|=i-1$. This generally identifies a unique
representation of the sequences as $A=X\seq{01}Y$ and $B=\zero X'\one Y'$
with $X\dom X'$ and $Y\dom Y'$. To reconstruct the unique tree $T$ mapped
to $(A,B)$, proceed recursively to identify $L$ from $X,X'$ and and $R$
from $Y,Y'$.

It remains to see that the cases where the above procedure fails to locate
a representation as above correspond exactly to cases where $T$ has only a
left or only a right sub-tree. One possibility is that $a(0)=0<1=b(1)$, and
then the above would give $|L|=-1$. In this case $A=\one Y$ and $B=\one Y'$
with $Y\dom Y'$, so the unique tree mapped to $A,B$ has only a right
sub-tree $R$, where $R$ is mapped to $Y,Y'$.

The other extreme is the case that $a(i)\ge b(i+1)$ for all $i<n$,
suggesting $|R|=-1$. In this case we have $A=X\zero$ and $B=\zero X'$ with
$X\dom X'$, and we find a tree with only a left sub-tree. Thus in all cases
the algorithm proceeds recursively to find a unique tree that is mapped to
$A,B$.
\end{proof}

\begin{lemma} \label{lem:w_id}
For a binary sequence $A$,
\[
W(A) = W(X) \chi_{A=X\zero} + W(Y) \chi_{A=\one Y}
     + \sum_{X\seq{01}Y=A} W(X)W(Y).
\]
\end{lemma}

Thus if $A$ ends with a \zero\ and equals $X\zero$, the RHS gets a
contribution of $W(X)$ from the first term. Similarly, if $A$ begins with a
\one\ there is a contribution from the second term. The sum in the RHS has
a term in the sum for each representation of $A$ as $X\seq{01}Y$. For
example, if $A=\seq{1011010}$ then we get
\[
\begin{aligned}
W(A) & = W(\seq{101101}) + W(\seq{011010})
                         + W(\one)W(\seq{1010}) + W(\seq{1011})W(\zero) \\
     & = 7 + 9 + 1\cdot5 + 2\cdot1
       = 23.
\end{aligned}
\]

\begin{proof}
By \corref{LR_w} it suffices to show that the RHS equals to the number of
binary trees that encode $A$. This is done by induction. If $A=X\zero$ ends
with a \zero, then there are $W(X)$ trees encoding $A$ with only a left
sub-tree (which can be any tree encoding $X$). If $A=\one Y$ then similarly
there are $W(Y)$ trees encoding $A$ with only a right sub-tree. Finally,
for any occurrence of \seq{01} in $A$ where $A=X\seq{01}Y$ there are $W(X)$
possible left sub-trees and $W(Y)$ possible right sub-trees, giving a term
in the sum.
\end{proof}

\section{Proof of \thmref{cycle}} \label{sec:cycle}
%%%%%%%%%%%%%%%%%%%%%%%%%%%%%%%%%

Consider a state $x$ of the exclusion process on the cycle. How many pairs
of sets $S,T$ collapse to $x$? Since the collapsing procedure begins by
placing the anti-particles at $T$, the unique $T$ is given by the set of
positions of anti-particles in $x$. There may be a number of different sets
$S$ that (together with $T$) collapse to the state $x$. In order for the
collapsing process to reach $x$ it is necessary that $S$ contains none of
the empty positions of $x$ (positions marked with \star's). The empty
positions in $x$ break up the cycle into a number of segments each
containing a sequence of particles (\one's) and anti-particles (\zero's).
Denote the binary segments of $x$ by $A_1,\ldots,A_l$.

During the collapsing procedure, if an element $p\in S$ results in a
particle being placed in some position $q$ to the left of $p$, there can be
no empty position in the interval $[q,p]$, since otherwise the particle
would have been placed there instead. Thus the elements of $S$ in each such
binary segment must collapse into the positions marked for particles in
that segment. It follows that for each binary segment $A_i$, the sequence
having \one's at the elements of $S$ in that segment is dominated by $A_i$,
and so there are $W(A_i)$ possibilities for the intersection of $S$ with
that segment. The total number of possibilities for $S$ is therefore $\prod
W(A_i)$, and the collapsed uniform measure of the state $x$ is
\[
\P(x) = \frac{\prod W(A_i)}{{N\choose a}{N\choose b}}.
\]

For example the cyclic state \seq{*10**10100*0101} may be reached from
\[
W(\seq{10})W(\phi)W(\seq{10100})W(\seq{0101}) = 2\cdot 1\cdot 9\cdot 2 = 36
\]
sets $S$ and so its probability is $36/{15\choose 5}{15\choose 6}$.

To show that the collapsed uniform measure is stationary, place at each
state $x$ a mass $m(x)=\prod W(A_i)$ --- a multiple of the collapsed
uniform measure. Let the mass flow according to the transition kernel of
the process (so if the process passes from $x$ to $y$ at rate $r$, mass
flows from $x$ to $y$ at a rate of $r\cdot m(x)$. It suffices to show that
the derivative of the mass at any state $x$ is 0.

Denote $x\to_e y$ if an action (sorting particles) along an edge $e$ leads
from state $x$ to state $y$. Since mass flowing from $x$ to itself makes no
difference, we only use this notation for $x\not = y$. Since the edge is
determined uniquely by $x$ and $y$, the edge subscript will usually be
omitted. Each edge is used at rate 1, so the mass derivative is given by
\begin{equation}\label{eq:deriv}
\frac{d}{dt}m(x) = \sum_{y\to_e x} m(y) - \sum_{x\to_e z} m(x)
                 = m(x) \sum_{y\to_e x} \left( \frac{m(y)}{m(x)}
                                             - \sum_{x\to_e z} 1 \right).
\end{equation}

We now associate each term in each of the sums of \eqref{eq:deriv} with one
of the binary sequences $A_i$ appearing in $x$ in such a way that the terms
associated with each sequence will cancel out, proving the Theorem. Since
edges connecting two empty positions lead from $x$ to itself and have been
disregarded, each edge intersects exactly one of the binary segments of
$x$, and the corresponding term is associated with that sequence.

The mass in the first sum --- flowing into $x$ --- corresponds to edges
with end-points marked '\seq{01}', '\seq{*1}', or '\seq{0*}' in $x$. For
each such edge we need to calculate the mass at the state resulting from
'unsorting' the edge. When such an edge is unsorted, the resulting state
$y$ is very similar to $x$. Indeed, $\{B_i\}$ are the binary sequences
appearing in $y$, then all but at one or two of them are equal to those in
$x$.

\begin{itemize}
\item
Consider first the case where $y\to_e x$ and the endpoints of $e$ are
marked '\seq{0*}' in $x$ and '\seq{*0}' in $y$. In this case $y$ has the
same binary sequences as $x$ except for two: $A_i = B_i\zero$ and $B_{i+1}
= \zero A_{i+1}$. Since $W(\zero A) = W(A)$, it follows that
$\frac{m(y)}{m(x)} = \frac{W(B_i)}{W(A_i)}$, where $B_i$ is $A_i$ with a
terminating \zero\ removed.
\item
Similarly, if $y\to_e x$ and the endpoints of $e$ are marked '\seq{*1}' in
$x$ and '\seq{1*}' in $y$, In this case $A_i = \one B_i$ and $B_{i-1} =
A_{i-1}\one$. Since $W(A\one) = W(A)$, we find that $\frac{m(y)}{m(x)} =
\frac{W(B_i)}{W(A_i)}$, where $B_i$ is $A_i$ with an initial \one\ removed.
\item
Finally, if $e$ is marked with '\seq{01}' in $x$, then spliting around that
edge we have $A_i=X\seq{01}Y$ and $B_i = X\seq{10}Y$, with all other
sequences beeing equal. Again, $\frac{m(y)}{m(x)} = \frac{W(B_i)}{W(A_i)}$.
\end{itemize}

Let us extend the $\to$ notation to binary sequences, so that $A\to B$ if
it is possible to pass from $A$ to $B$ by either removing an initial \zero,
removing a terminating \one, or replacing a \seq{10} by \seq{01} somewhere
in $A$. Consider the terms in the RHS of \eqref{eq:deriv} that are
associated with the binary sequence $A_i$ of $x$. After substituting the
above for $\frac{m(y)}{m(x)}$, these terms come to
\[
m(x) \left( \sum_{B\to A_i} \frac{W(B)}{W(A_i)} - \sum_{A_i\to C} 1 \right),
\]
and so it suffices to prove for an arbitrary sequence $A$ that
\begin{equation}\label{eq:w_id}
\sum_{B\to A} W(B) = \sum_{A\to C} W(A).
\end{equation}

Given a binary sequence $A$, associated terms in the RHS correspond to each
occurrence of '\seq{10}' in the sequence, as well as to an initial \zero\
if there is one and to terminating \one\ if the sequence ends with a \one.
The number of such terms is always one more than the number of times
'\seq{01}' appears in $A$. 

The terms in the LHS, are $W(B)$ where $B$ results from $A$ either by
replacing '\seq{01}' by '\seq{10}' at some place, or by removing an initial
\one\ or a terminating \zero\ if $A$ has them. Consider a pair of sequences
$B\to A$ where $A=X\seq{01}Y$ and $B=X\seq{10}Y$. Since $B\dom A$, $B$ also
dominates any sequence that $A$ dominates. The difference $W(B)-W(A)$
corresponds to sequences that $B$ dominates and $A$ does not. Such a
sequence must be of the form $X'\seq{10}Y'$, where $X\dom X'$ and $Y\dom
Y'$, and so $W(B)=W(A)+W(X)W(Y)$.

Substituting this for of $W(B)$ in \eqref{eq:w_id} results in cancellation
of all but one term on the RHS, and the resulting needed identity is
exactly that given by \lemref{w_id}.

\section{The infinite setting} \label{sec:line}
%%%%%%%%%%%%%%%%%%%%%%%%%%%%%%

\begin{proof}[Proof of \thmref{line}]
It is known (see \cite{FKS}) that there is a unique non-trivial stationary
measure with marginals $p,q$ for particles and anti-particles, and that
those are all the non-trivial extremal stationary measures. It only
needs to be shown that the collapsed i.i.d. measures are stationary  
and have the correct marginals.

The collapsed i.i.d. measure is the limit as $N\to\infty$ of the collapsed
uniform measure on a cycle of length $N$ with $[pN]$ particles and $[qN]$
anti-particles. Since in the finite case particles and anti-particles
are uniformly 
distributed over all subsets of the appropriate size, in the
limit they have densities $p$ and $q$ respectively.
Since correlations in this
measure decay exponentially in the distance, it follows that the limit is
stationary.
\end{proof}

The collapsing procedure sheds new light some known results.
First, the fact that second class particles induce a factoring of the
stationary distribution (see \cite{DJLS}): given that 0 is an empty position, the state in
$\Z^+$ and in $\Z^-$ are independent. With the collapsing procedure,
the state in $\Z^+$ depends only on the positive elements of $S,T$. 
Conditioned on the event that 0 remains empty, no particle crosses
from $\Z^+$ to $\Z^-$, and so the state on $\Z^-$ is determined by the
negative elements of $S,T$ and is therefore independent of the state
on $\Z^+$.

Next, consider the relation demonstrated in \cite{speer,FFK} to a certain biased random
walk on $\Z$. For marginals $p$ for particles
and $q$ for anti-particles, the corresponding random walk has steps of
$-1,0,1$ with step distribution given by $\P(X=1)=(1-p)(1-q)$, and
$\P(X=-1)=pq$. When $p+q<1$ we have $\E X = 1-p-q>0$.

The random walk is very naturally coupled with the pair $S,T$ of i.i.d.\
subsets of $\Z$ with densities $p,q$, and is given by $Z_n =
|T^c\cap[-n,0)| - |S\cap[-n,0)|$ (note that the sets are explored backwards), which clearly has the above step
distribution. Lemma~2.5 of \cite{FFK} says that the distribution of the
distance between second class particles is the same as the hitting time of
1 by the random walk. In our notation this is a statement regarding the
distance between empty positions. Given $S,T$ such that there is no
particle at 0, the next hole to the left of 0 is at $-n$ exactly when the
corresponding $Z$ hits $1$ at time $n$. 

The following is standard and easily seen: The stationary measures
seen from a single second class particles are 
derived from the stationary measures with anti-particles, by the
following rule. Condition on having an empty position at 0,
and place a second class particle there. All particles are first class
particles, and any empty spot to the right of 0 is also filled with a
first class particle. All other positions become empty. By selecting
the densities of particles and anti-particles the asymptotic densities
to either side can be controlled.

Extend the definition of the random walk to negative indices by
$Z_{-n} = - |T^c\cap[0,n)| + |S\cap[0,n)|$. Theorem~5 of \cite{FFK}
states that the stationary measure for the process seen from a single
second class particle with given asymptotic densities is very close to
a product of independent i.i.d.\ measures on $\Z^+$ and $\Z^-$ with
the given densities. Formally, the stationary measure and the product
measure  can be coupled 
so that they are exponentially unlikely to differ in many positions. 

This again follows from the collapsing procedure, since the final
location of a particle is exponentially unlikely to be far from the
position given by $S$.

\section{Open problems} \label{sec:open}
%%%%%%%%%%%%%%%%%%%%%%%

The following conjectures are a generalization of our results to more
classes of particles, both in qualitative and quantitative terms.

\begin{conj}
In a cycle with particles of a number of classes, in the stationary
distribution, the least likely states are those where if the cycle is cut
at some point the particles are arranged in reversed order of speed.
\end{conj}

Let $x$ be one of the above states (there are $N$ of them).

\begin{conj}
In the stationary distribution, $\P(x)=\prod {N\choose s_i}^{-1}$, where
$n$ is the cycle length and $s_i$ is the number of particles of class at
least $i$.
\end{conj}

\begin{conj}
In the stationary distribution, the probability of any other state is an
integral multiple of $\P(x)$.
\end{conj}

Finally, is there a useful generalization of the collapsing procedure for
processes with more classes of particles? The obvious generalizations of
repeated collapsing do not appear to give the correct stationary
distribution. The smallest case where they fail is that of a cycle of
length 4 with all different particles (or 3 particles and an empty space).
It turns out that $\mu(1324) \not = \mu(1423)$ under the stationary measure
$\mu$.


\begin{thebibliography}{1}

\bibitem{BBHM}
I.~Benjamini, N.~Berger, C.~Hoffman, and E.~Mossel.
\newblock Mixing times of the biased card shuffling and the asymmetric
  exclusion process.
\newblock To appear in Trans. of AMS., 2004.

\bibitem{DJLS}
B.~Derrida, S.~Janowsky, J.~Lebowitz, and E.~Speer.
\newblock Exact solution of the totally asymmetric simple exclusion process:
  shock profiles.
\newblock {\em J. Stat. Phys.}, 73:833--874, 1993.

\bibitem{DS}
E.~Duchi and G.~Schaeffer.
\newblock A combinatorial approach to jumping particles.
\newblock To appear, J. Comb. The. A., 2004.

\bibitem{FFK}
P.~Ferrari, L.~Fontes, and Y.~Kohayakawa.
\newblock Invariant measures for a two species asymmetric process.
\newblock {\em J. Stat. Phys.}, 76(5):1153--1177, 1994.

\bibitem{FKS}
P.~Ferrari, C.~Kipnis, and E.~Saada.
\newblock Microscopic structure of traveling waves in the asymmetric simple
  exclusion process.
\newblock {\em Ann. Prob.}, 19:226--244, 1991.

\bibitem{liggett1}
T.~Liggett.
\newblock {\em Interacting Particle Systems}.
\newblock Fundamental Principles of Mathematical Sciences, 276.
  Springer-Verlag, New York, 1985.

\bibitem{liggett2}
T.~Liggett.
\newblock {\em Stochastic Interacting Systems: Contact, Voter and Exclusion
  Processes}.
\newblock Fundamental Principles of Mathematical Sciences, 324.
  Springer-Verlag, New York, 1999.

\bibitem{speer}
Speer.
\newblock The two species asymmetric simple exclusion process.
\newblock In A.~V. M.~Fannes, C.~Maes, editor, {\em Micro, Meso, and
  Macroscopic Approaches in Physics, Procedings of the NATO AR Workshop "On
  Three Levels"}, Leuven, July 1993 1994.

\end{thebibliography}
\end{document}